\documentclass[11pt,a4paper]{amsart}

\usepackage{amssymb,amsmath,amsthm}

\newtheorem{theorem}{Theorem}[]
\newtheorem{corollary}{Corollary}

\newtheorem{claim}{Claim}

\usepackage{kbordermatrix} 
\renewcommand{\kbldelim}{(} 
\renewcommand{\kbrdelim}{)}

\usepackage{color,latexsym,amsmath,amssymb,path, amsthm,enumitem}
\usepackage{graphicx}
\usepackage{amssymb}
\usepackage{enumerate}

\renewcommand{\geq}{\geqslant}
\renewcommand{\leq}{\leqslant}

\newcommand{\etc}{,\dots,}
\newcommand{\C}{\mathbb{C}}
\newcommand{\R}{\mathbb{R}}

\newcommand{\la}{\lambda}
\newcommand{\al}{\alpha}
\newcommand{\be}{\beta}
\newcommand{\si}{\sigma}
\newcommand{\ga}{\gamma}

\begin{document}

\title{Gershgorin disks for multiple eigenvalues of non-negative matrices}

\dedicatory{Dedicated to the memory of Ji\v{r}\'i Matou\v{s}ek}

\author{Imre B\'ar\'any}
\address{Alfr\'ed R\'enyi Mathematical Institute of Mathematics,
Hungarian Academy of Sciences
P.O. Box 127, 1364 Budapest,  Hungary, and
 Department of Mathematics,
 University College London,
Gower Street, London WC1E 6BT, UK}

\author{J\'ozsef Solymosi}
\address{Department of Mathematics, University of  British Columbia,  1984 Mathematics Road,
Vancouver, BC, Canada V6T 1Z2}

\maketitle

\begin{abstract}
Gershgorin's famous circle theorem states that all eigenvalues of a square matrix lie in disks (called Gershgorin disks) around the diagonal elements. Here we show that if the matrix entries are non-negative and an eigenvalue has geometric multiplicity at least two, then this eigenvalue lies in a smaller disk. The proof uses geometric rearrangement inequalities on sums of higher dimensional real vectors which is another new result of this paper.
\end{abstract}

\section{Introduction and main result}
Gershgorin's circle theorem \cite{GE} is a fundamental and widely used result on localizing the eigenvalues of square matrices. It states that all eigenvalues are in disks (called Gershgorin disks) around the diagonal elements.

The main goal of this paper is to improve Gershgorin's theorem under special conditions, namely, when the matrix is non-negative and has a multiple eigenvalue. We show that such an eigenvalue lies in disks of smaller radius around a diagonal element. For the proof we establish various geometric inequalities concerning rearrangements of vector sums. This is an interesting connection between convex geometry and matrix theory. The geometric point of view in eigenvalue problems is certainly not new but this particular connection seems to be new.

 Here we show that if the matrix entries are non-negative and an eigenvalue has geometric multiplicity at least two, then this eigenvalue lies in a smaller disk.

Let $D(a,r)$ denote the disk with center $a$ and radius $r$ on the complex plane:

\[
D(a,r)=\left\{x \in \C : |x-a|\leq r\right\}.
\]

For an $n\times n$ complex matrix, $A=[a_{ij}],$ the Gershgorin disks are $D(a_{ii},R_i)$ where $R_i=\sum_{j: i\neq j}|a_{ij}|.$ The most commonly cited form of Gershgorin's theorem says that every eigenvalue of $A$ lies in some $D(a_{ii},R_i)$.  Varga's nice book {\em Gershgorin and His Circles} \cite{VA} surveys various applications and extensions of this important theorem. An interesting and recent theorem of Marsli and Hall \cite{MH1} states that
if an eigenvalue of a matrix $A$ has geometric multiplicity $k,$ then it lies in at least $k$ of the Gershgorin disks of $A.$ They have extended this result in subsequent papers \cite{FMH, MH2, MH3, MH4}. Here we focus on the $k=2$ case for non-negative matrices.

Understanding the spectra of a matrix is a central question both in applied and pure mathematics. Here 
are some facts and results. There are two particular eigenvalues for which the multiplicity is of great importance; the largest eigenvalue which determines the spectral radius of the matrix and the multiplicity of the eigenvalue ``0'' since it determines the rank of the matrix. There are also applications using the smallest eigenvalue. For example Roy shows in \cite{RO} that the Euclidean representation number of a graph is closely related to the multiplicity of the smallest eigenvalue. The multiplicity of the largest and the second largest eigenvalues play a key role in some numerical methods.  Del Corso \cite{DC} considers the problem of approximating an eigenvector belonging to the largest eigenvalue by the so called power method. It is proved that the rate of convergence depends on the ratio of the two largest eigenvalues and on their multiplicities. The rate increases with the multiplicity of the largest eigenvalue and decreases with the multiplicity of the second eigenvalue. In graph theory the Colin de Verdi\`ere number is the multiplicity of the second largest eigenvalue of the adjacency matrix, maximized by weighting the edges and nodes. For more details and the exact definition we refer to the papers \cite{LO1} and \cite{HLS}.

Gershgorin's circle theorem is intertwined with the Perron-Frobenius theory. It is one of the tools used to bound the spectral radius of a matrix. It follows from the Perron-Frobenius theorem that the largest magnitude eigenvalue of any non-negative matrix is a positive real number, see in e.g. ~\cite{BP}.

Let us define the {\em half Gershgorin} disks, $D(a_{ii},r_i),$ which are subsets of the original. Instead of $ R_i=\sum_{j: i\neq j}|a_{ij}|$ we take the partial sum of the $\lfloor n/2\rfloor$ largest terms. This sum is denoted by $r_i$.

Recall that the {\sl geometric multiplicity} of an eigenvalue $\lambda$ of $A$ is the dimension of the corresponding eigenspace of $A$, that is, the kernel of $A-\lambda I$. (Its algebraic multiplicity is the multiplicity of the root $\lambda$ of the polynomial $\det (A-xI)$.)

We are going to show that multiple geometric eigenvalues are in smaller Gershgorin disks when the matrix is non-negative.

\begin{theorem}\label{th:halfG}
Let $A=\{a_{ij}\}$ be an $n\times n$ non-negative (real) matrix and $\lambda$ an eigenvalue of $A$ with geometric multiplicity at least two. Then $\lambda$ is in a half Gershgorin disk, $D(a_{ii},r_i),$ for some $i.$
\end{theorem}

Actually we are going to prove that such an eigenvalue lies in the disk $D(a_{ii},r)$ and various values of $r$ for some suitable $i$. The proofs are based on geometric estimates that are of independent interest. They are given in the next section.

\section{Rearrangement inequalities for vectors}

Assume $V=\{v_1\etc v_n\} \subset \R^d$  and $\sum_1^nv_i=0$. Further, let $\al_1 \ge \ldots \ge \al_n \ge 0$ be real numbers. We write $[n]$ for the set $\{1\etc n\}$.

\begin{theorem}\label{th:first} Under the above conditions set $\beta= \al_{\lfloor n/2 \rfloor+1}$. Then for every permutation $\si$ of $[n]$
\[
\| \sum _1^{n} \al_i v_{\si(i)}\| \le \max_{i\in [n]}\|v_i\| \sum_1^n |\al_i-\beta|.
\]
\end{theorem}

\begin{corollary}\label{cor:first}Under the above conditions, for every permutation $\si$ of $[n]$
\[
\| \sum _1^{n} \al_i v_{\si(i)}\| \le \max_{i\in [n]}\|v_i\| \sum_1^{\lfloor n/2 \rfloor}\al_i.
\]
\end{corollary}

In the second geometric estimate we need a technical assumption.

\begin{theorem}\label{th:second} Let $V=\{v_1,\ldots,v_n\}\subset \R^d$ satisfy the previous assumption. Suppose further that the $v_i$ are ordered with decreasing (Euclidean) length,
that is, $\|v_1\|\ge \ldots \ge \|v_{n}\|$. Let $\ga\in [\al_{j+1},\al_j]$ for some $j \in [n-1]$. Then for every permutation $\si$ of $[n]$
\[
\| \sum _1^{n} \al_i v_{\si(i)}\| \le \sum_1^j \al_i\|v_i\|- \frac {\ga}2\left[\sum_1^j\|v_i\|-\sum_{j+1}^{n}\|v_i\|\right].
\]
\end{theorem}

Here of course one wants to choose $j$ and $\ga$ so that the right hand side is as small as possible. When $j=\lceil n/2 \rceil$, the sum between the brackets is non-negative. Choosing any $\ga$ from the interval $[\al_{j+1},\al_j]$ gives the following.

\begin{corollary}\label{cor:second}Under the above conditions for every permutation $\si$ of $[n]$
\[
\| \sum _1^{n} \al_i v_{\si(i)}\| \le \sum_1^{\lceil n/2 \rceil} \al_i\|v_i\|.
\]
\end{corollary}

We mention that the estimates in Theorems \ref{th:first} and \ref{th:second} are incomparable; sometimes the first, other times the second gives the better bound.

\section{Proof of the rearrangement inequalities}

{\bf Proof} of Theorem~\ref{th:first}. First fix some $\gamma \ge 0$. Then
\[
\sum_1^n\al_iv_{\si(i)}=\sum_1^n\al_iv_{\si(i)}-\sum_1^n\ga v_{\si(i)}=\sum_1^n(\al_i-\ga)v_{\si(i)}.
\]
By the triangle inequality
\[
\|\sum_1^n\al_iv_{\si(i)}\|\le \max_{i\in [n]}\|v_i\| \sum_1^n|\al_i-\gamma|.
\]
Set $k=\lfloor n/2 \rfloor$ and define $\be=\al_{k+1}$. It can be proven that the function $\gamma \to \sum_1^n|\al_i-\gamma|$ takes its minimum at $\gamma=\beta$ when $n$ is odd, and at every $\ga$ from the interval $[\al_{k+1},\al_k]$ when $n$ is even. $\Box$

\medskip
Corollary~\ref{cor:first} follows immediately since with the above $k$ and $\be$
\begin{eqnarray*}
\sum_1^n|\al_i-\be|&=&\sum_1^k(\al_i-\be)+\sum_{k+1}^n(\be -\al_i)\\
&=& \sum_1^k\al_i -\sum_{\ell}^n\al_i \le \sum_1^k\al_i
\end{eqnarray*}
where $\ell$ equals $k+1$ for even $n$ and $k+2$ for odd $n$.

\medskip
{\bf Proof} of Theorem~\ref{th:second}. The zonotope $Z(V)$ spanned by $V$ is, by definition, the set
\[
Z(V)=\left\{\sum_{i\in [n]}\xi_iv_i: 0\le \xi_i \le 1\; (\forall i)\right\}.
\]
Let $B$ denote the Euclidean unit ball of $\R^d$. We claim first that
\begin{equation}\label{eq:zono}
Z(V) \subset \frac 12 \Big(\|v_1\|+\cdots +\|v_{n}\|\Big)B.
\end{equation}
It is well-known \cite{McM} and easy to check that $Z(V)$ is the convex hull of the points $s(W)=\sum_{v\in W}v$ where $W \subset V$.
Thus it suffices to show that for every $W \subset V$, $\|s(W)\|\le  \frac 12 (\|v_1\|+\ldots +\|v_{n}\|)$. Fix $U \subset V$ such that $s(U)$ has maximal length among all $s(W)$. Set $z=s(U)$ and observe that $-z=s(V \setminus U)$ as $s(V)=0$. Since $\|z\|=\|-z\|$ evidently, we have
\[
2\|z\|=\|z\|+\|-z\|=\|s(U)\|+\|s(V\setminus U)\| \le \sum_1^n\|v_i\|
\]
by the triangle inequality. This implies that $\|z\|\le \frac 12 \sum_1^n\|v_i\|$.

\medskip
We observe next that
\begin{eqnarray*}
\sum _1^n \al_i v_{\si(i)}&=&\sum_1^j(\al_i-\ga)v_{\si(i)}+ \sum_1^j\ga v_{\si(i)}+ \sum_{j+1}^{n}\al_iv_{\si(i)}\\
    &=&\sum_1^j(\al_i-\ga)v_{\si(i)}+ \ga \left[\sum_1^jv_{\si(i)} +\sum_{j+1}^{n}\frac {\al_i}{\ga}v_{\si(i)}\right].
\end{eqnarray*}

The expression between the brackets is a vector $u$ in $Z(V)$ so $\|u\| \le \frac 12 \sum_1^{n}\|v_i\|$. By the triangle inequality the norm of $\sum_1^{n} \al_i v_{\si(i)}$ is at most
\begin{eqnarray*}
\sum_1^j(\al_i-\ga)\|v_{\si(i)}\|+ \ga \|u\|&\le& \sum_1^j(\al_i-\ga)\|v_i\|+  \frac {\ga}2 \sum_1^{n}\|v_i\|\\
   &=& \sum_1^j\al_i\|v_i\| -\frac {\ga}2 \left[\sum_1^j\|v_i\|-\sum_{j+1}^{n}\|v_i\|\right].
\end{eqnarray*}
\qed

\section{Proof of Theorem~\ref{th:halfG}}

We first recall the simple proof of Gershgorin's original theorem. Let $v=(v_1,\ldots,v_n)$ be an eigenvector with eigenvalue $\lambda$ where $v_i$ are complex numbers. Assume $|v_i|=\max_{j\in [n]}|v_j|$. Then $\sum_{j=1}^na_{ij}v_j=\lambda v_i$ implying
\begin{equation}\label{eq:basic}
(\lambda-a_{ii})v_i=\sum_{j:j\ne i}a_{ij}v_j.
\end{equation}
Taking absolute value on both sides and using $|v_i|\ge |v_j|$ shows that $\lambda \in D(a_{ii},R_i)$ with $R_i=\sum_{j:j\ne i}a_{ij}$ indeed.

When the eigenvalue $\lambda$ has geometric multiplicity at least two, then its eigenspace contains a nonzero vector $v=(v_1,\ldots,v_n)$ whose components sum to zero: $\sum_1^nv_i=0$. Indeed, let $u$ and $w$ be two linearly independent eigenvectors from the eigenspace of $\lambda$. If $\sum_1^nu_i=0$, then $v=u$ is a suitable eigenvector. If not, then $v=(\sum_1^nw_i)u-(\sum_1^nu_i)w$ has the required property.

As any multiplier of $v$ is still an eigenvector, we can suppose that the largest magnitude component of $v,$ $v_i$, is a positive real number. Actually we can and do assume that $v_i=1$. Then the other components, $v_j$, are complex numbers with $|v_j|\le 1$.

\medskip
The {\bf proof} of Theorem~\ref{th:halfG} is based on equation (\ref{eq:basic}) plus the condition that $\sum_1^nv_j=0$. As $\C$ is a vector space of dimension 2 over $\R$, we can consider the components $v_j$ of $v$ as vectors in $\R^2$. Then Theorem~\ref{th:first} with $d=2$ applies to the $v_j \in \R^2$, we just have to imagine that on the right hand side of (\ref{eq:basic}) $v_i$ is added with coefficient zero. So define $b_{ii}=0$ and $b_{ij}=a_{ij}$ if $i\ne j$. Let $b^*$ be the median of the sequence $b_{i1},\ldots,b_{in}$. Theorem~\ref{th:first} gives then that $\lambda$ lies in the disk $D(a_{ii},r)$ where
\begin{equation}\label{eq:strong}
r=\sum_{j\ne i}|b_{ij}-b^*|.
\end{equation}

\smallskip
The proof of Theorem~\ref{th:halfG} uses Corollary~\ref{cor:first}: $\lambda$ lies in the disk $D(a_{ii},r)$ where $r$ is the sum of the largest $\lfloor n/2 \rfloor$ entries in the $i$th row of $A$ (disregarding $a_{ii})$.  Note that in general the estimate in (\ref{eq:strong}) is gives a better bound on $r$ than Theorem~\ref{th:halfG}. $\Box$

\smallskip
We can also apply Corollary \ref{cor:second} to the components of $v$, considered again as vectors in $\R^2$. This gives that $\la$ lies in the disk $D(a_{ii},r)$ where $r$ is the sum of the $k=\lceil n/2 \rceil$ largest entries in row $i$ of $A$ (disregarding $a_{ii}$ again). In any special case a better estimate may come from the more general Theorem~\ref{th:second}.

\medskip
{\bf Remark 1.}
One could hope that an eigenvalue with (geometric) multiplicity 3 or higher should lie strictly inside the half Gershgorin disk. The simple example below shows that this is not the case.

Let $A$ be an $n\times n$ matrix with $n=3k$, consisting of three $k\times k$ blocks along the main diagonal, with each block being a doubly stochastic matrix. Then $\lambda=1$ is an eigenvalue with multiplicity 3, which lies on the boundary of each half Gershgorin disk $D(a_{ii},r_i)$. Indeed $r_i$ is the sum of the largest $\lfloor n/2 \rfloor$ entries of the $i$th row (disregarding $a_{ii}$) which equals $1-a_{ii}$.

This example shows, however, that $\lambda$ lies in the ``third Gershgorin disk''. This is the disk centred at $a_{ii}$ and of radius $r$ which is the sum of the largest $n/3$ entries in the $i$th row (disregarding again $a_{ii}$). We return to this question at the end of the paper.

\section{Examples}

In what follows we show examples illustrating the limits of possible extensions of the results above. Note that one can not expect in general that a multiple eigenvalue is strictly inside the half Gershgorin disk. The simp\-lest illustration to this is the matrix $A$ below where $1$ is an eigenvalue with (geometric) multiplicity two.

\[A=\left[
\begin{array}{ccc}
 0 & 1 & 1  \\
 1 & 0 & 1  \\
 1 & 1 & 0  \\
\end{array}
\right]\]

Next we are going to give further examples. The first two show that Theorem~\ref{th:halfG} does not extend to real matrices that have both positive and negative entries. The second is a positive semidefinite Hermitian matrix (with complex entries) where the triple eigenvalue ``0'' lies on the boundary of the half Gershgorin disk. Perhaps some form of Theorem~\ref{th:halfG} can be extended to such matrices.

\subsection{Real matrices with both positive and negative entries} The matrices in Theorem~\ref{th:halfG} have non-negative entries. This condition cannot be deleted as the following symmetric circulant matrix with $0,\pm 1$ entries shows:
\[
B = \left[\begin{array}{ccccc}0&1&-1&-1&1\\1&0&1&-1&-1\\-1&1&0&1&-1\\-1&-1&1&0&1\\1&-1&-1&1&0\end{array} \right]
\]
 Like every $5\times 5$ symmetric circulant matrix, $B$ has two multiple eigenvalues. They are $\sqrt{5}\approx 2.236$ and $-\sqrt{5}$ and both lie outside the half Gershgorin disk.

 The following $7 \times 7$ matrix is again circulant and has $0,\pm 1$ entries. Its multiple eigenvalue $\approx -3.494$ is even further from the half Gershgorin disk which has radius 3 around the origin.

 \[C=\left[
\begin{array}{ccccccc}
 0 & 1 & -1 & 1 & 1 & -1 & 1 \\
 1 & 0 & 1 & -1 & 1 & 1 & -1 \\
 -1 & 1 & 0 & 1 & -1 & 1 & 1 \\
 1 & -1 & 1 & 0 & 1 & -1 & 1 \\
 1 & 1 & -1 & 1 & 0 & 1 & -1 \\
 -1 & 1 & 1 & -1 & 1 & 0 & 1 \\
 1 & -1 & 1 & 1 & -1 & 1 & 0 \\
\end{array}
\right]\]

\subsection{A positive semidefinite matrix}
The next construction gives a $9\times 9$ positive semidefinite Hermitian matrix $H$ with the triple eigenvalue ``0'' lying on the boundary of the half Gershgorin disk. (This is very different from the example in Remark 1 where the half disk and the third disk were the same.) The other eigenvalue is 6 and it lies on the ``quarter disk''. This example comes from the Hesse configuration of 9 points and 12 lines in ${\mathbb{CP}}^2$ \cite{HE}. The matrix $H$ looks interesting on its own right. It shows further that strengthening Theorem \ref{th:halfG} to more general matrices (with high multiplicity eigenvalues) might be difficult.

One possible realization of the Hesse configuration is given by the following 9  points on the complex projective plane
\[
\begin{array}{ccc}
p_1 = (0, 1, -1) & p_2 = (0, 1, -\omega) & p_3 = (0, 1, -\omega^2)\\
p_4 = (1, 0, -1) & p_5 = (1, 0, -\omega^2) & p_6 = (1, 0, -\omega)\\
p_7 = (1, -1, 0) & p_8 = (1, -\omega, 0) & p_9 = (1, -\omega^2, 0)\\
\end{array}
\]
where $\omega=\frac{-1+i\sqrt{3}}{2}$ is a third root of unity. In this arrangement each point lies
on four lines and each line contains three points. Our first matrix, $A,$ records the linear dependencies of the points. It has 9 columns, one for each point,  and 12 rows, one for each line. If $p_i,p_j$ and $p_k$ are collinear, then there are nonzero complex multipliers $\alpha,\beta,\gamma$ such that $\alpha p_i+\beta p_j+\gamma p_k=0.$ For example the sixth (highlighted) row in the matrix $A$ below represents the equation
\[-\omega^2(0, 1, -1)-(0, 1, -\omega)-\omega (0, 1, -\omega^2)=(0,0,0).
\]

Thus the matrix $A$ encodes the linear dependencies of collinear triples in the point-line arrangement of the Hesse configuration.

\[
A=\left[
\begin {array}{ccccccccc} 1&0&0&-1& 0&0&1& 0&0\\
\noalign{\medskip}0&0&1& 0&-1& 0&1& 0& 0\\
\noalign{\medskip}0&1&0& 0& 0&-1&1& 0& 0\\
\noalign{\medskip} 0& 0& 0& 0& 0& 0&-\omega^2&-1&-\omega\\
\noalign{\medskip} 0& 0& 0&-\omega^2&-\omega&-1& 0& 0& 0\\
\noalign{\medskip}{\bf -\boldsymbol{\omega}^2}&{\bf -1}&{\bf -\boldsymbol{\omega}}& {\bf 0}& {\bf 0}& {\bf 0}& {\bf 0}& {\bf 0}& {\bf 0}\\
\noalign{\medskip} 0&\omega& 0& 0&-1& 0& 0&1& 0\\
\noalign{\medskip} 0& 0&-\omega^2& 0& 0&1& 0& 0&-1\\
\noalign{\medskip}-\omega^2& 0& 0& 0&1& 0& 0& 0&-1\\
\noalign{\medskip}\omega& 0& 0& 0& 0&-1& 0&1& 0\\
\noalign{\medskip} 0&1& 0&-\omega& 0& 0& 0& 0&\omega\\
\noalign{\medskip} 0& 0&\omega&-1& 0& 0& 0&1& 0\end {array}\right]
\]

The points of the Hesse configuration satisfy the homogeneous system of equations  $A\bf{x}=0$ where $x_i\in {\mathbb{CP}}^2$. An affine image of a solution is also a solution, implying that the rank of $A$ is at most 6. It is easy to see that the rank is exactly 6: the rank remains the same if one multiplies a matrix with its Hermitian transpose (complex conjugate transpose). So consider the $9\times9$ matrix $H=\overline{A^T}A$.

\[ H=\left[ \begin {array}{ccccccccc}  \phantom{-}4&\phantom{-}\omega\phantom{^2}&\phantom{-}\omega^2&-1 &  -\omega\phantom{^2}& -\omega^2& \phantom{-}1 &\phantom{-}\omega^2&\phantom{-}\omega\phantom{^2}\\ \noalign{\medskip}\phantom{-}\omega^2& \phantom{-}4&\phantom{-}\omega\phantom{^2}&  -\omega\phantom{^2}& -\omega^2&-1 & \phantom{-}1 &\phantom{-}\omega^2&\phantom{-}\omega\phantom{^2}\\ \noalign{\medskip}\phantom{-}\omega\phantom{^2}&\phantom{-}\omega^2& \phantom{-}4&  -\omega^2&-1 &  -\omega\phantom{^2}& \phantom{-}1 &\phantom{-}\omega^2&\phantom{-}\omega\phantom{^2}\\ \noalign{\medskip}-1 &  -\omega^2&  -\omega\phantom{^2}& \phantom{-}4&\phantom{-}\omega^2&\phantom{-}\omega\phantom{^2}&-1 &-1 &-1 \\ \noalign{\medskip}  -\omega^2&  -\omega\phantom{^2}&-1 &\phantom{-}\omega\phantom{^2}& \phantom{-}4 &\phantom{-}\omega^2&-1 &-1 &-1 \\ \noalign{\medskip}  -\omega\phantom{^2}&-1 &  -\omega^2&\phantom{-}\omega^2&\phantom{-}\omega\phantom{^2}& \phantom{-}4&-1 &-1 &-1 \\ \noalign{\medskip} \phantom{-}1 & \phantom{-}1 & \phantom{-}1 &-1 &-1 &-1 & \phantom{-}4&\phantom{-}\omega\phantom{^2}&\phantom{-}\omega^2\\ \noalign{\medskip}\phantom{-}\omega\phantom{^2}&\phantom{-}\omega\phantom{^2}&\phantom{-}\omega\phantom{^2}&-1 &-1 &-1 &\phantom{-}\omega^2& \phantom{-}4&\phantom{-}\omega\phantom{^2}\\ \noalign{\medskip}\phantom{-}\omega^2&\phantom{-}\omega^2&\phantom{-}\omega^2&-1 &-1 &-1 &\phantom{-}\omega\phantom{^2}&\phantom{-}\omega^2& \phantom{-}4\end {array} \right]
 \]

Matrix $H$ is a positive semidefinite Hermitian matrix that has two eigenvalues: 0 with multiplicity 3 (so the rank of $A$ is indeed 6) and 6 with multiplicity 6. All non-diagonal entries have norm one and the diagonal entries are 4. Thus $\lambda=0$ is on the boundary of the half Gershgorin disk $D(4,4)$ and $\lambda=6$ on the boundary of $D(4,2)$, the ``quarter disk''.

\begin{figure} [h]
    \centering
\includegraphics[ width=10cm]{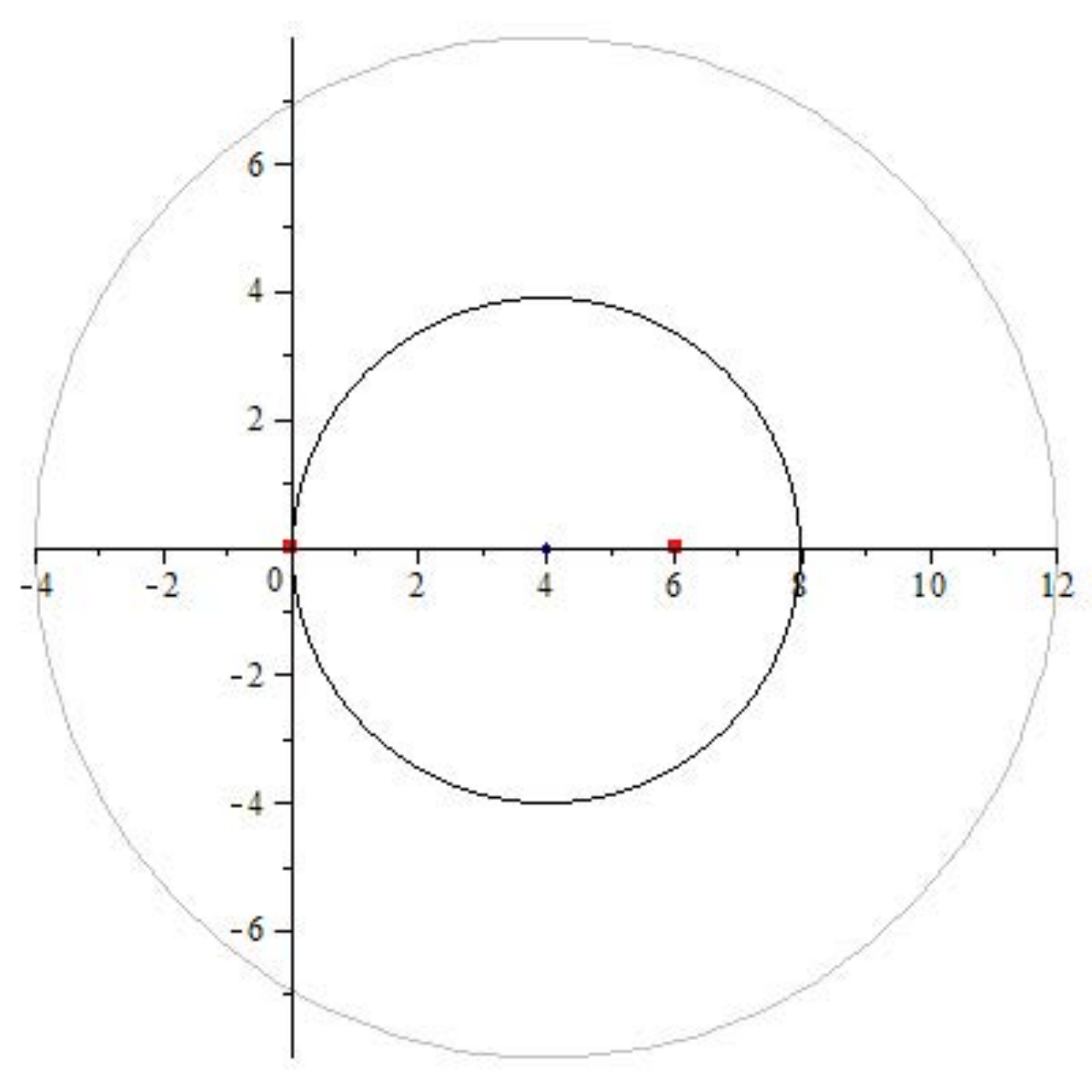}
\caption{The Gershgorin disk and half-disk of $H$}
    \label{fig:Hesse}
\end{figure}

\section{Remarks}
There are several questions that remain open.
\begin{itemize}
\item What can be said about the location of an eigenvalue with larger multiplicity? Our method, using the zonotope $Z(V)$ in the proof of Theorem~\ref{th:second} has its limitations. Perhaps inequality (\ref{eq:zono}) can be improved.  For instance, for an eigenvalue with multiplicity at least $k$ one would like to use an eigenvector $v=(v_1,\ldots,v_n)$ such that the corresponding zonotope $Z(V)$ satisfies
\[
Z(V) \subset c \left(\|v_1\|+\ldots +\|v_{n}\|\right)B
\]
where $c$ decreases as $k$ grows. Unfortunately one can not expect $c$ to go below $\frac1 \pi$, see Exercise 14.9 in \cite{St})
\item How about other matrices? What is the radius of the shrunken Gershgorin disk which contains  a multiple eigenvalue of a general complex matrix? Are there better bounds for special matrices, like real or positive semidefinite Hermitian matrices?
\end{itemize}

\section*{Acknowledgment}
This research was supported by ERC Advanced Research Grant no 267165 (DISCONV). Imre B\'ar\'any is partially supported by Hungarian National Research Grant K 111827. J\'ozsef Solymosi is partially supported by Hungarian National Research Grant NK 104183 and by an NSERC Discovery Grant. We are indebted to three anonymous referees for very useful comments and information that have improved the presentation of this paper.

\end{document}